\documentclass[12pt,a4paper,twoside]{article}

\usepackage{amsmath,amsfonts,amsthm,amsopn,color,amssymb}
\usepackage{palatino}
\usepackage{graphicx}

\newcommand{\eps}{\varepsilon}

\newcommand{\R}{\mathbb{R}}

\newcommand{\RN}{{\mathbb{R}^N}}
\newcommand{\RT}{{\mathbb{R}^3}}

\DeclareMathOperator{\supp}{supp}

\renewcommand{\le}{\leqslant}
\renewcommand{\ge}{\geqslant}
\renewcommand{\a }{\alpha }

\renewcommand{\b }{\beta }

\renewcommand{\d }{\delta }

\newcommand{\g }{\gamma }

\newcommand{\n }{\nabla }

\renewcommand{\t}{\theta}
\renewcommand{\O}{\Omega}

\newcommand{\G}{\Gamma}

\newcommand{\Ne}{\mathcal{N}}
\newcommand{\E}{\mathcal{E}}

\newcommand{\M}{\mathcal{M}}

\newcommand{\D }{{\mathcal D}^{1,2}(\RT)}
\renewcommand{\H}{H^{1}(\RT)}
\newcommand{\irn }{\int_{\RN}}
\newcommand{\irt }{\int_{\RT}}

\def\bbm[#1]{\mbox{\boldmath $#1$}}

\newtheorem{theorem}{Theorem}[section]
\newtheorem{lemma}[theorem]{Lemma}

\newtheorem{remark}[theorem]{Remark}

\renewenvironment{proof}{\noindent{\textbf{Proof\quad}}}{$\hfill\square$\vspace{0.2 cm}\\}
\newenvironment{proofmain}{\noindent{\textbf{Proof of Theorem  \ref{main}\quad}}}{$\hfill\square$\vspace{0.2 cm}\\}

\newenvironment{proofmain3}{\noindent{\textbf{Proof of Theorem  \ref{main3}\quad}}}{$\hfill\square$\vspace{0.2 cm}\\}
\newenvironment{proofmain4}{\noindent{\textbf{Proof of Theorem  \ref{main4}\quad}}}{$\hfill\square$\vspace{0.2 cm}\\}
\newenvironment{proofnonex}{\noindent{\textbf{Proof of Theorem \ref{nonex}\quad}}}{$\hfill\square$\vspace{0.2 cm}\\}

\title{{\bf Ground state solutions for the nonlinear Schr\"odinger-Maxwell
equations}}

\author{A. Azzollini \thanks{Dipartimento di Matematica, Universit\`a degli
Studi di Bari,  Via E. Orabona 4, I-70125 Bari, Italy, e-mail: {\tt azzollini@dm.uniba.it}}
 \; \& \;
A. Pomponio\thanks{Dipartimento di Matematica, Politecnico di Bari, Via Amendola 126/B, I-70126 Bari, Italy, e-mail:
{\tt a.pomponio@poliba.it}}}

\date{}

\begin{document}

\maketitle
\begin{abstract}
In this paper we study the nonlinear Schr\"odinger-Maxwell equations
\[
\left\{
\begin{array}{ll}
-\Delta u + V(x)u + \phi u = |u|^{p-1}u  & \hbox{in }\RT,
\\
-\Delta \phi = u^2 & \hbox{in }\RT.
\end{array}
\right.
\]
If $V$ is a positive constant, we prove the existence of a ground state solution $(u,\phi)$ for $2<p<5$. The non-constant potential case is treated under suitable geometrical assumptions on $V$, for $3<p<5$. Existence and non-existence results are proved also when the nonlinearity exhibits a critical growth.
\end{abstract}

\tableofcontents

\section{Introduction}

In this paper we consider the problem
\begin{equation}\label{ScMa}\tag{$\mathcal{SM}$}
\left\{
\begin{array}{ll}
-\Delta u + V(x)u + \phi u = f'(u)  & \hbox{in }\RT,
\\
-\Delta \phi = u^2 & \hbox{in }\RT,
\end{array}
\right.
\end{equation}
where $V:\RT \to \R$ and $f\in C^1(\RT,\R).$ Such a system
represents the nonlinear Schr\"odinger-Maxwell equations in the
electrostatic case. In \cite{BF}, the potential $V$ has been
supposed constant, and the linear version of the problem (i.e. $f
\equiv 0$) has been studied as an eigenvalue problem for a bounded domain.
The linear Schr\"odinger-Maxwell equations have been treated also in \cite{C1,CG}, where the potential $V$ has been supposed radial.

The nonlinear case has been considered in \cite{AR,C2,DM,DA,Ru},
where existence and multiplicity results  have been stated when
$V$ is a positive constant. By means of the Pohozaev's fibering method, a
multiplicity result has been proved in \cite{S} also in the
non-homogeneous case, that is when a non-homogeneous term
$g(x)\in L^2(\RT)$ is added on the right hand side of the first
equation of \eqref{ScMa} (see also \cite{CS}). On the other hand,
nonexistence results for \eqref{ScMa} can be found in
\cite{DM2,Ru}. For a related problem see \cite{PS}.

In this paper we will look for {\it ground state solutions} to the
problem \eqref{ScMa}, namely for couples $(u,\phi)$ which solve
\eqref{ScMa} and minimize the action functional associated to
\eqref{ScMa} among all possible solutions. The problem of finding
such a type of solutions is a very classical problem: it has been
introduced by Coleman, Glazer and Martin in \cite{CGM}, and
reconsidered by Berestycki and Lions in \cite{BL1} for a class of
nonlinear equations including the Schr\"odinger's one. Later on
the existence and the profile of ground state solutions have been
studied for a plethora of problems by many authors; of course we
can not mention all these results.

In the first part of the paper, we are interested in considering pure power type
nonlinearities so that the problem we will deal with becomes
\begin{equation}\label{SM}
\left\{
\begin{array}{ll}
-\Delta u + V(x)u + \phi u = |u|^{p-1}u  & \hbox{in }\RT,
\\
-\Delta \phi = u^2 & \hbox{in }\RT,
\end{array}
\right.
\end{equation}
where $2<p<5$. The solutions $(u,\phi)\in \H \times \D$ of
\eqref{SM} are the critical points of the action functional $\mathcal{E}
\colon \H \times \D \to \R$, defined as
\[
\mathcal{E}(u,\phi):=\frac 12 \irt |\n u|^2 + V(x)u^2
-\frac 14 \irt |\n \phi|^2
+\frac 12 \irt \phi u^2
-\frac 1{p+1} \irt |u|^{p+1}.
\]
We are interested in finding a ground state solution of \eqref{SM}, that is a solution $(u_0,\phi_0)$ of \eqref{SM} with the property of having the least action among all possible solutions of \eqref{SM}, namely $\E(u_0,\phi_0)\le \E(u,\phi)$, for any solution $(u,\phi)$ of \eqref{SM}.

The action functional $\E$ exhibits a strong indefiniteness, namely it is
unbounded both from below and from above on infinite dimensional
subspaces. This indefiniteness can be removed using the reduction
method described in \cite{BFMP}, by which we are led to study a
one variable functional that does not present such a strongly
indefinite nature.

The main difficulty related with the problem of finding
the critical points of the new functional, consists in the lack of
compactness of the Sobolev spaces embeddings in the unbounded
domain $\RT$. Usually, at least when $V$ is radially symmetric, such a difficulty is overcome by
restricting the functional to the natural constraint of the radial
functions where compact embeddings hold. In particular, in
\cite{DM} a radial solution having minimal energy among all the
radial solutions has been found. However we are not able to say if
that solution actually is a ground state for our equation. This is
the reason why we will use an alternative method, based on a
concentration-compactness argument on suitable measures, to
recover compactness.

We analyze two different situations. First we assume that $V$ is a positive constant and, following an idea of
Ruiz \cite{Ru}, we look for a minimizer of the reduced functional
restricted to a suitable manifold $\M$. Such a manifold has two
interesting features: it is a natural constraint
for the reduced functional and it contains, in a
sense that we will explain later (see Remark \ref{re:grst}), every solution of the problem
\eqref{SM}. The main result we get is the
following
\begin{theorem}\label{main}
If $V$ is a positive constant, then the problem \eqref{SM} has a ground state
solution for any $p\in ]2,5[.$
\end{theorem}

In fact, it is standard to see that such a ground state solutions does not change sign, so we can assume it positive.

Then we study \eqref{SM} assuming the following hypotheses on $V:$
\begin{itemize}
\item[({\bf V1})] $V\in C(\RT,\R)$; \item[({\bf V2})] $0<C_1\le
V(x)\le C_2$, for all $x\in \RT$; \item[({\bf V3})]
$V_\infty:=\liminf_{|y|\to\infty}V(y) \ge V(x)$, for all $x\in
\RT$, and the inequality is strict for some $x\in \RT$.
\end{itemize}
These kind of hypotheses on the potential were introduced by
Rabinowitz \cite{R} to study the nonlinear Schr\"odinger equation
\[
-\Delta u + V(x)u  = f'(u) \qquad \hbox{in }\RT.
\]
Because of technical difficulties we are not allowed to use the same device as in the constant potential case. We study the reduced functional restricted to the Nehari manifold and we are able to prove the existence result only for $3<p<5$:
\begin{theorem}\label{main2}
If $V$ satisfies ({\bf V1-3}) then the problem \eqref{SM} has a
ground state solution for any $p\in ]3,5[$.
\end{theorem}

Theorems \ref{main} and \ref{main2} will be proved in Section \ref{se:subcr}.

In the second part of the paper we consider the critical case, namely the
case when the nonlinearity presents at infinity the same
behavior of the power $t^{2^*-1}$, where $2^*=6$ is the critical
exponent for the Sobolev embeddings in dimension $3$. Here a further obstacle to compactness
arises, in fact, it is well known that the embedding of the space
$H^{1}(\O)$ into the Lebesgue space $L^{2^*}(\O)$ is
not compact, even if $\O$ is a bounded set in $\RT.$

The problem becomes
\begin{equation}\label{SMcr}
\left\{
\begin{array}{ll}
-\Delta u + V(x)u + \phi u = u^5  & \hbox{in }\RT,
\\
-\Delta \phi = u^2 & \hbox{in }\RT.
\end{array}
\right.
\end{equation}

By \cite{DM2}, we have the following
\begin{theorem}[D'Aprile \& Mugnai \cite{DM2}]\label{th:DM2}
Suppose that $V$ is a positive constant.
Let $(u,\phi)\in \H\times \D$ be a solution of the problem \eqref{SMcr}, then $u=\phi=0$.
\end{theorem}

We extend this nonexistence result to the case of a non-constant potential $V$. We prove the following nonexistence theorem, based on a Pohozaev-type identity.
\begin{theorem}\label{nonex}
Suppose that $V$ satisfies ({\bf V2}) and
\begin{itemize}
\item[({\bf V4})] $V\in C^1(\RT,\R)$; \item[({\bf V5})] $0\le (\n
V(x)\mid x) \le C_3$, for all $x\in \RT$.
\end{itemize}
Let $(u,\phi)\in \H\times \D$ be a solution of the problem \eqref{SMcr}, then $u=\phi=0$.
\end{theorem}

Then, in the same spirit of \cite{BN} (see also \cite{C} for the
Klein-Gordon-Maxwell equation), we add a lower order perturbation
to the first equation of \eqref{SMcr}, namely we look for
solutions to the system
\begin{equation}\label{SMc}
\left\{
\begin{array}{ll}
-\Delta u + V(x)u + \phi u = |u|^{q-1}u+u^5 & \hbox{in }\RT,
\\
-\Delta \phi = u^2 & \hbox{in }\RT,
\end{array}
\right.
\end{equation}
where $q\in]3,5[$. The solutions $(u,\phi)\in \H \times \D$ of \eqref{SMc}
are the critical points of the action functional $\mathcal{E}^* \colon \H \times \D \to \R$, defined as
\begin{align*}
\mathcal{E}^*(u,\phi):=&\;\frac 12 \irt |\n u|^2 + V(x)u^2
-\frac 14 \irt |\n \phi|^2
+\frac 12 \irt \phi u^2
\\
&-\frac 1{q+1} \irt |u|^{q+1}
-\frac 16 \irt u^6.
\end{align*}
The effect of the additive perturbation is to lower the energy.
This causes that the ground state level of the functional falls
into an interval where compactness holds. As a consequence we get
the following two results, respectively for the constant and the non-constant potential case:

\begin{theorem}\label{main3}
Let $V$ be a positive constant. Then the problem \eqref{SMc} has a ground state
solution.
\end{theorem}

\begin{theorem}\label{main4}
Let $V$ satisfy ({\bf V1-3}). Then the problem \eqref{SMc} has a ground state
solution.
\end{theorem}

We will prove these three last theorems in Section \ref{se:crit}.

\vspace{0.5cm}
\begin{center}
{\bf NOTATION}
\end{center}

\begin{itemize}
\item For any $1\le s< +\infty$, $L^s(\RT)$ is the usual Lebesgue space endowed with the norm
\[
\|u\|_s^s:=\irt |u|^s;
\]
\item $\H$ is the usual Sobolev space endowed with the norm
\[
\|u\|^2:=\irt |\n u|^2+ u^2;
\]
\item $\D$ is completion of $C_0^\infty(\RT)$ with respect to the norm
\[
\|u\|_{\D}^2:=\irt |\n u|^2;
\]
\item for any $r>0,$ $x\in\RT$ and $A\subset \RT$
\begin{align*}
B_r(x) &:=\{y\in\RT\mid |y-x|\le r\},\\
B_r    &:=\{y\in\RT\mid |y|\le r\},\\
A^c    &:= \RT\setminus A;
\end{align*}
\item $C,\,C',\,C_i$ are positive constants which can change from
line to line;
\item $o_n(1)$ is a quantity which goes to zero as $n \to +\infty$.
\end{itemize}

\section{The subcritical case}\label{se:subcr}

\subsection{Some preliminary results}\label{se:prelim}

We first recall some well-known facts (see, for instance
\cite{BF,C1,C2,CG,DM,Ru}). For every $u\in L^{12/5}(\RT)$, there
exists a unique $\phi_u\in \D$ solution of
\[
-\Delta \phi=u^2,\qquad \hbox{in }\RT.
\]
It can be proved that $(u,\phi)\in H^1(\RT)\times \D$ is a
solution of \eqref{SM} if and only if $u\in\H$ is a critical point
of the functional $I\colon \H\to \R$ defined as
\begin{equation}\label{eq:defI}
I(u)= \frac 12 \irt |\n u|^2 + V(x) u^2 +\frac 14 \irt \phi_u u^2
-\frac{1}{p+1}\irt |u|^{p+1},
\end{equation}
and $\phi=\phi_u$.\\
The functions $\phi_u$ possess the following properties (see
\cite{DM} and \cite{Ru})
\begin{lemma}\label{le:prop}
For any $u\in\H$, we have:
\begin{itemize}
\item[i)] $\|\phi_u\|_{\D}\le C \|u\|^2,$ where $C$ does not
depend from $u.$ As a consequence there exists $C'>0$ such that
$$
\irt\phi_u u^2\le C'\|u\|_\frac {12}5^4;
$$
\item[ii)] $\phi_u\ge 0;$
\item[iii)] for any $t>0$:
$\phi_{tu}=t^2\phi_u;$
    \item[iv)] for any $\t>0$: $\phi_{u_\t}(x)=\t^2\phi_u(\t x)$,
    where $u_\t(x)=\t^2 u(\t x)$;
    \item[v)] for any $\O\subset\RT$
measurable,
$$
\int_{\O}\phi_u u^2=\int_{\O}\irt\frac{u^2(x)u^2(y)}{|x-y|}dx\,dy.
$$
\end{itemize}
\end{lemma}

\subsection{The constant potential case}

In this section we will assume that $V$ is a positive constant. Without lost of generality, we suppose $V\equiv 1.$ It can be proved (see \cite{DM2,Ru}) that
if $(u,\phi)\in \H\times \D$ is a solution of \eqref{SM}, then it
satisfies the following Pohozaev type identity
\begin{equation}\label{eq:Poho}
\irt \frac12 |\n u|^2 + \frac 32 u^2 + \frac 54 \phi u^2 - \frac 3{p+1}|u|^{p+1}=0.
\end{equation}
As in \cite{Ru}, we introduce the following manifold
\[
\M:=\left\{u\in \H\setminus \{0\} \;\Big{|}\; G(u)=0 \right\},
\]
where
$$
G(u):=
\irt \frac 32 |\n u|^2 + \frac 12 u^2 +
\frac 34 \phi_u u^2 - \frac {2p-1}{p+1}|u|^{p+1}.
$$
\begin{remark}\label{re:grst}
Observe that if $u\in \H$ is a nontrivial critical point of $I$,
then $u\in \M$, since $G(u)=0$ can be obtained by a linear
combination of $\langle I'(u), u\rangle =0$ and \eqref{eq:Poho},
with $\phi=\phi_u$. As a consequence if $(u,\phi)\in \H\times\D$ is a solution of \eqref{SM}, then $u\in \M$.
\end{remark}

The next lemma describes some properties of the manifold $\M$:
\begin{lemma}\label{le:M}
\begin{enumerate}
\item For any $u\in \H$, $u\neq 0$, there exists a unique number $\bar \t>0$
such that $u_{\bar \t}\in \M.$ Moreover
\[
I(u_{\bar \t})=\max_{\t \ge 0} I(u_\t);
\]
\item there exists a positive constant $C$, such that for all $u
\in \M$, $\|u\|_{p+1}\ge C$;
\item $\M$ is a natural constraint of $I$, namely every critical point of $I|_{\M}$ is a critical point for $I.$
\end{enumerate}
\end{lemma}
\begin{proof}
We refer to \cite{Ru}. In particular, as regards point 3, we have
to point out that Ruiz \cite{Ru} has just proved that the minimum
of $I|_{\M}$ is in fact a critical point of $I$: the same
arguments can be adapted to prove that $\M$ is a natural
constraint of $I$.
\end{proof}

By 3 of Lemma \ref{le:M} we are allowed to look for critical
points of $I$ restricted to $\M$.
\\
Moreover, by 1 of Lemma \ref{le:M}, the map
$\t:\H\setminus\{0\}\to\R_+$ such that for any $u\in\H,$ $u\neq
0:$
\begin{equation*}
I\left( u_{\t (u)}\right)=\max_{\t\ge 0} I( u_\t)
\end{equation*}
is well defined.
\\
Set
\begin{align}
c_1 &= \inf_{g\in \G} \max_{\t\in [0,1]} I(g(\t)); \nonumber
\\
c_2 &= \inf_{u\neq 0} \max_{\t \ge 0} I( u_\t);\nonumber
\\
c_3 &= \inf_{u\in\M} I(u);\nonumber
\end{align}
where
\begin{equation}\label{eq:gamma}
\G=\left\{g\in C\big([0,1],\H\big) \mid g(0)=0,\;I(g(1))\le 0,
\;g(1)\neq 0\right\}.
\end{equation}
\begin{lemma}\label{le:ccc}
The following equalities hold
\[
c:=c_1=c_2=c_3.
\]
\end{lemma}
\begin{proof}
Taking into account 1 of Lemma \ref{le:M} and the fact that
for small $\|u\|$ we have (see \cite[Theorem 3.2, Step1]{Ru})
\begin{equation*}
\irt \frac 32 |\n u|^2 + \frac 12 u^2 +
\frac 34 \phi_u u^2
> \irt \frac {2p-1}{p+1}|u|^{p+1},
\end{equation*}
the conclusion follows using the same arguments of
\cite[Proposition 3.11]{R}.
\end{proof}

\begin{remark}\label{re:grst2}
By point 3 of Lemma \ref{le:M} and Remark \ref{re:grst}, we argue that if $u\in \M$ is such that $I(u)=c$, then $(u,\phi_u)$ is a ground state solution of \eqref{SM}.
\end{remark}

\subsubsection{Proof of Theorem \ref{main}}

Let $(u_n)_n \subset \M$ such that
\begin{equation}\label{eq:lim}
\lim_n I(u_n)=c.
\end{equation}
We define the functional $J\colon \H\to\R$ as:
\begin{equation*}
J(u)=\irt \frac {p-2}{2p-1} |\n u|^2 + \frac{p-1}{2p-1}  u^2 +
\frac{p-2}{2(2p-1)} \phi_{u} u^2.
\end{equation*}
Observe that for any $u\in\M,$ by $ii$ of Lemma \ref{le:prop} we
have $I(u)=J(u)\ge 0.$
\\
By \eqref{eq:lim}, we deduce that $(u_n)_n$ is bounded in $\H,$ so
there exists $\bar u\in\H$ such that, up to a subsequence,
\begin{align}
&u_n\rightharpoonup \bar u\quad\hbox{weakly in }\H,    \label{eq:weakM}
\\
&u_n\to \bar u\quad\hbox{in }L^s(B), \hbox{with }B\subset \RT, \hbox{bounded, and }1\le s<6.  \nonumber 
\end{align}
To prove Theorem \ref{main}, we need some compactness on the
sequence $(u_n)_n.$ To this end, we use a
concentration-compactness argument on the positive measures so
defined: for every $u_n \in \H$,
\begin{equation}\label{eq:meas}
\nu_n(\O)= \int_\O \frac {p-2}{2p-1}|\n u_n|^2
+ \frac{p-1}{2p-1} u_n^2
+  \frac{p-2}{2(2p-1)} \phi_{u_n} u_n^2.
\end{equation}
By \eqref{eq:lim} we have
\[
\nu_n(\RT)=J(u_n)\to c
\]
and then, by P.L.~Lions \cite{L1}, there are three possibilities:
\begin{description}
\item[\;\;{\it vanishing}\,{\rm :}] for all $r>0$
\[
\lim_n \sup_{\xi \in \RT}\int_{B_r(\xi)} d \nu_n =0;
\]
\item[\;\;{\it dichotomy}\,{\rm :}] there exist a constant $\tilde
c\in (0, c)$, two sequences $(\xi_n)_n$ and $(r_n)_n$, with $r_n
\to +\infty$ and two nonnegative measures $\nu_n^1$ and $\nu_n^2$
such that
\begin{align*}
0\le \nu_n^1 + \nu_n^2 \le \nu_n,&\qquad \nu_n^1(\RT) \to \tilde
c,\;\; \nu_n^2(\RT) \to c -\tilde c,
\\
\hbox{supp}(\nu_n^1)\subset B_{r_n}(\xi_n),&\qquad
\hbox{supp}(\nu_n^2)\subset \RT \setminus B_{2r_n}(\xi_n);
\end{align*}
\item[\;\;{\it compactness}\,{\rm :}] there exists a sequence
$(\xi_n)_n$ in $\RT$ with the following property: for any $\d>0$,
there exists $r=r(\d)>0$ such that
\[
\int_{B_r(\xi_n)} d \nu_n \ge c -\d.
\]
\end{description}
Arguing as in \cite{WZ}, we prove the following
\begin{lemma}\label{le:concentr}
Compactness holds
for the sequence of measures $(\nu_n)_n$, defined in \eqref{eq:meas}.
\end{lemma}
\begin{proof}
{\sc Vanishing does not occur}
\\
Suppose by contradiction, that for all $r>0$
\[
\lim_n \sup_{\xi \in \RT}\int_{B_r(\xi)} d \nu_n =0.
\]
In particular, we deduce that there exists $\bar r>0$ such that
\begin{equation*}
\lim_n \sup_{\xi \in \RT}\int_{B_{\bar r}(\xi)} u_n^2=0.
\end{equation*}
By \cite[Lemma I.1]{L2}, we have that $u_n\to 0$ in $L^s(\RT),$ for $2< s <6.$ As a consequence, since $(u_n)_n \subset \M$ and by Lemma \ref{le:prop}, we get
\begin{align*}
0\le I(u_n)
&\le \irt \frac 32 |\n u_n|^2 + \frac 12 u_n^2
+\frac 14 \phi_{u_n} u_n^2
-\frac{1}{p+1}|u_n|^{p+1}
\\
&=-\frac 12 \irt \phi_{u_n}u_n^2
+ \frac {2p-2}{p+1}  \irt |u_n|^{p+1} \to0
\end{align*}
which contradicts \eqref{eq:lim}.
\\
\\
{\sc Dichotomy does not occur}
\\
Suppose by contradiction that there exist a constant $\tilde c\in
(0, c)$, two sequences $(\xi_n)_n$ and $(r_n)_n$, with $r_n \to
+\infty$ and two nonnegative measures $\nu_n^1$ and $\nu_n^2$ such
that
\begin{align*}
0\le \nu_n^1 + \nu_n^2 \le \nu_n,&\qquad \nu_n^1(\RT) \to \tilde
c,\;\; \nu_n^2(\RT) \to c -\tilde c,
\\
\hbox{supp}(\nu_n^1)\subset B_{r_n}(\xi_n),&\qquad
\hbox{supp}(\nu_n^2)\subset \RT \setminus B_{2r_n}(\xi_n).
\end{align*}
Let $\rho_n \in C^1(\RT)$ be such that $\rho_n\equiv 1$ in
$B_{r_n}(\xi_n)$, $\rho_n\equiv 0$ in $\RT \setminus
B_{2r_n}(\xi_n)$, $0\le \rho_n\le 1$ and $|\n \rho_n|\le 2/r_n$.
\\
We set
\[
v_n:=\rho_n u_n,\qquad w_n:=(1-\rho_n)u_n.
\]
It is easy to see that
\begin{align*}
\liminf_n J(v_n) &\ge \tilde c,
\\
\liminf_n J(w_n) &\ge c - \tilde c.
\end{align*}
Moreover, denoting $\O_n:= B_{2r_n}(\xi_n)\setminus
B_{r_n}(\xi_n)$, we have
\begin{equation*}
\nu_n(\O_n) \to 0, \qquad \hbox{as }n\to \infty,
\end{equation*}
namely
\begin{align}
\int_{\O_n}|\n u_n|^2 +  u^2_n \to 0, \qquad \hbox{as }n\to \infty,\nonumber 
\\
\int_{\O_n}\phi_{u_n} u_n^2 \to0, \qquad \hbox{as }n\to \infty.\label{eq:2}
\end{align}
By simple computations, we infer also
\begin{align}
\int_{\O_n}|\n v_n|^2 +  v^2_n \to 0, \qquad \hbox{as }n\to \infty,\nonumber
\\
\int_{\O_n}|\n w_n|^2 +  w^2_n \to 0, \qquad \hbox{as }n\to \infty.\nonumber
\end{align}
Hence, we deduce that
\begin{align}
\irt |\n u_n|^2 +  u^2_n
&=\irt|\n v_n|^2 +  v^2_n
\quad+\irt|\n w_n|^2 +  w^2_n
+o_n(1),\label{eq:norm}
\\
\irt |u_n|^{p+1}
&=\irt|v_n|^{p+1}+\irt|w_n|^{p+1}
+o_n(1).\label{eq:p+1}
\end{align}
Moreover, by point $v$ of  Lemma \ref{le:prop} and \eqref{eq:2}, we have
\begin{align}
\irt \phi_{u_n} u_n^2 &= \irt \phi_{v_n} v_n^2 + \irt \phi_{w_n} w_n^2
+ 2\int_{B_{r_n}}\!\int_{B_{2r_n}^c}\!\!\!\!\frac{u_n^2(x)u_n^2(y)}{|x-y|}d x\,d y + o_n(1) \nonumber
\\
&\ge \irt \phi_{v_n} v_n^2 + \irt \phi_{w_n} w_n^2 + o_n(1). \label{eq:phi}
\end{align}
Hence, by \eqref{eq:norm} and \eqref{eq:phi}, we get
\begin{align}
J(u_n) \ge J(v_n)+J(w_n)+o_n(1).\nonumber
\end{align}
Then
\begin{align*}
c=\lim_n J(u_n)
\ge\liminf_n J(v_n)
+\liminf_n J(w_n)
\ge \tilde c+(c-\tilde c)=c,
\end{align*}
hence
\begin{align}
\lim_n J(v_n)&=\tilde c, \label{eq:jv}
\\
\lim_n J(w_n)&=c-\tilde c. \nonumber
\end{align}
We recall the definition of the functional $G\colon \H\to \R$
\[
G(u)=\irt \frac 32|\n u|^2+\frac 12 u^2+\frac 3 4\phi_u u^2- \frac{2p-1}{p+1}|u|^{p+1}
\]
and that if $u\in \M$, then $G(u)=0$. By \eqref{eq:norm}, \eqref{eq:p+1}  and \eqref{eq:phi}, we have
\begin{equation}
0=G(u_n) \ge G(v_n)+G(w_n)+o_n(1)\label{eq:G}.
\end{equation}
By Lemma \ref{le:M}, for any $n \ge 1$, there exists $\t_n>0$ such that $ ({v_n})_{\t_n} \in \M$, and then
\begin{equation}\label{eq:nehari}
\irt\frac{3}2 \t_n^2|\n v_n|^2 + \frac 12 v_n^2 + \frac{3} 4
\t_n^2 \phi_{v_n}v_n^2 = \irt\frac{2p-1}{p+1}\t_n^{2p-2}
|v_n|^{p+1}.
\end{equation}
We have to distinguish three cases.
\\
\
\\
{\sc Case 1:} up to a subsequence, $G(v_n) \le 0$.
\\
By \eqref{eq:nehari} we have
\begin{align*}
\irt\frac 32 (\t_n^{2p-2}-\t_n^2) |\n v_n|^2 + \frac 12 (\t_n^{2p-2}-1) v_n^2 + \frac 3 4(\t_n^{2p-2}-\t_n^2)\phi_{v_n}v_n^2 \le 0,
\end{align*}
which implies that $\t_n\le 1$. Therefore, for all $n\ge 1$
\begin{align*}
c \le I\big( ({v_n})_{\t_n}\big) = J\big( ({v_n})_{\t_n}\big)
\le J(v_n) \to \tilde c < c,
\end{align*}
which is a contradiction.
\\
\
\\
{\sc Case 2:} up to a subsequence, $G(w_n) \le 0$.
\\
We can argue as in the previous case.
\\
\
\\
{\sc Case 3:} up to a subsequence, $G(v_n) > 0$ and $G(w_n) > 0$.
\\
By \eqref{eq:G}, we infer that $G(v_n)=o_n(1)$ and $G(w_n)=o_n(1)$.
If $\t_n\le 1+o_n(1)$, we can repeat the arguments of Case 1.
Suppose that
\begin{equation*}
\lim_n \t_n=\t_0>1.
\end{equation*}
We have
\begin{align*}
o_n(1)&= G(v_n) = \irt \frac 32|\n v_n|^2+\frac 12 v_n^2+\frac 3 4\phi_{v_n} v_n^2- \frac{2p-1}{p+1}|v_n|^{p+1}
\\
&=\irt\frac 32\left( 1-\frac{1}{\t_n^{2p-4}}\right) |\n v_n|^2+\frac 12 \left(1-\frac{1}{\t_n^{2p-2}}\right)v_n^2\\
&\quad+\irt\frac 34\left( 1-\frac{1}{\t_n^{2p-4}}\right) \phi_{v_n} v_n^2
\end{align*}
and so $v_n \to 0$ in $\H$, but we get a contradiction with \eqref{eq:jv}.
\\
Hence we conclude that dichotomy can not occur.
\end{proof}
Now we are able to yield the following
\\

\begin{proofmain}
Let $(u_n)_n$ be a sequence in  $\M$ such that \eqref{eq:lim}
holds.
\\
We define the measures $(\nu_n)_n$ as in \eqref{eq:meas}; by Lemma \ref{le:concentr}
there exists a sequence $(\xi_n)_n$
in $\RN$ with the following property: for any $\d>0$, there exists
$r=r(\d)>0$ such that
\begin{equation}\label{eq:bc}
\int_{B_r^c(\xi_n)} \frac{p-2}{2p-1} |\n u_n|^2 + \frac{p-1}{2p-1}
                u_n^2
                +  \frac{p-2}{2(2p-1)} \phi_{u_n} u_n^2< \d.
\end{equation}
We define the new sequence of functions
$v_n:=u_n(\cdot-\xi_n)\in\H.$ It is easy to see that
$\phi_{v_n}=\phi_{u_n}(\cdot-\xi_n)$, and hence $v_n\in\M.$
Moreover, by \eqref{eq:bc}, we have that for any $\d>0$, there
exists $r=r(\d)>0$ such that
\begin{equation}\label{eq:small}
\|v_n\|_{H^1(B_r^c)}<\d\;\hbox{uniformly for $n\ge 1$}.
\end{equation}
Since, by \eqref{eq:weakM}, $(v_n)_n$ is bounded in $\H,$ certainly there
exist a subsequence (likewise labelled) and $\bar v\in\H$
such that
\begin{align}
&v_n\rightharpoonup \bar v\quad\hbox{weakly in }\H, \label{eq:weak2}
\\
&v_n\to \bar v \quad\hbox{in }L^s(B), \hbox{with }B\subset \RT, \hbox{bounded, and }1\le s<6.  \label{eq:loc2}
\end{align}
By \eqref{eq:small}, \eqref{eq:weak2} and \eqref{eq:loc2}, we have that, taken $s\in [2,6[$,
for any $\d>0$ there exists $r>0$ such that, for any $n\ge 1$
large enough
\begin{align*}
\|v_n-\bar v\|_{L^s(\RT)}
& \le \|v_n-\bar v\|_{L^s(B_r)}+\|v_n-\bar v\|_{L^s(B^c_r)}
\\
& \le \d + C\left(\|v_n\|_{H^{1}(B_r^c)}+\|\bar v\|_{H^{1}(B_r^c)}\right)\le
(1+2C)\d,
\end{align*}
where $C>0$ is the constant of the embedding
$H^{1}(B_r^c)\hookrightarrow L^s(B^c_r).$
We deduce that
\begin{equation}\label{eq:conver}
v_n\to \bar v\hbox{ in }L^s(\RT),\;\hbox{for any }s\in [2,6[.
\end{equation}
Since $\phi$ is continuous from $L^{12/5}(\RT)$ to $\D,$
from \eqref{eq:conver} we deduce that
\begin{align}
\phi_{v_n}  \to \phi_{\bar v} \;\hbox{ in }\D,& \qquad \hbox{as } n\to \infty,  \nonumber\\
\irt\phi_{v_n}v_n^2  \to  \irt\phi_{\bar v}\bar v^2,& \qquad
\hbox{as } n\to \infty.  \label{eq:convphi3}
\end{align}
Since $(v_n)_n$ is in $\M$, by 2 of Lemma \ref{le:M}
$(\|v_n\|_{p+1})_n$ is bounded below by a positive constant. As a
consequence, \eqref{eq:conver} implies that $\bar v\neq 0.$
Proceeding as in \cite[Theorem 3.2, Step 4]{Ru}, by \eqref{eq:conver} and
\eqref{eq:convphi3} we can show that $v_n\to\bar v$ in $\H$ so
that $\bar v\in\M$ and $I(\bar v)=c.$ By Remark \ref{re:grst2}, we have that $(\bar v,\phi_{\bar v})$ is a ground state
solution of \eqref{SM}.
\end{proofmain}

\subsection{The non-constant potential case}

In this section we suppose that the potential $V$ satisfies ({\bf
V1-3}) and that $p\in]3,5[.$

In order to get our result, we will use a very standard device: we
will look for a minimizer of the functional \eqref{eq:defI}
restricted to the Nehari manifold
\begin{equation*}
\Ne=\left\{u\in\H\setminus\{0\}\mid \widetilde
G(u)=0\right\},
\end{equation*}
where
\[
\widetilde G(u):= \irt  |\n u|^2 +  u^2 +  \phi_u u^2 -
|u|^{p+1}.
\]
The following lemma describes some properties of the Nehari
manifold~$\Ne$:
\begin{lemma}\label{le:N}
\begin{enumerate}
\item For any $u\neq 0$ there exists a unique number $\bar t>0$ such that $\bar t u\in \Ne$ and
\[
I(\bar t u)=\max_{ t \ge 0}I( t u);
\]
\item there exists a positive constant $C$, such that for all $u
\in \Ne$, $\|u\|_{p+1}\ge C$;
\item $\Ne$ is a $C^1$ manifold.
\end{enumerate}
\end{lemma}
\begin{proof}
Points 1 and 2 can be proved using standard arguments (see, for example, \cite{R}).
\\
3. Observe that for any $u\in\H$ we have
$$
\widetilde G(u)= 4 I(u) -\irt\left(|\n u|^2 + V(x) u^2\right) -\frac{p-3}{p+1} \irt
|u|^{p+1},
$$
and then, by point 2, for any $u\in\Ne$ we have
\begin{equation*}
\langle \widetilde G'(u), u\rangle=-2\irt\left(|\n u|^2 + V(x) u^2\right) - (p-3) \irt
|u|^{p+1} \le - C < 0.
\end{equation*}
\end{proof}

The Nehari manifold $\Ne$ is a natural constrained for the
functional $I,$ therefore we are allowed to look for critical
points of $I$ restricted to $\Ne$.

In view of this, we assume the following definition
\begin{equation*}
c_V:=\inf_{u\in\Ne} I(u),
\end{equation*}
so that our goal is to find $\bar u\in\Ne$ such that $I(\bar u)=c_V$, by which we would deduce that $(\bar u,\phi_{\bar u})$ is a ground state solution of \eqref{SM}.

First we recall some preliminary lemmas which can be obtained by
using the same arguments as in \cite{R} (see also \cite{AP}).

As a consequence of the Lemma \ref{le:N}, we are allowed to
define the map $t:\H\setminus\{0\}\to\R_+$ such that for any
$u\in\H,$ $u\neq 0:$
\begin{equation*}
I\left( t(u)u\right)=\max_{t\ge 0} I(t u).
\end{equation*}
\begin{lemma}\label{le:ccc2}
The following equalities hold
\begin{equation*}
c_V=\inf_{g\in \G} \max_{t\in [0,1]}
I(g(t))=\inf_{u\neq 0} \max_{t \ge 0} I( tu),
\end{equation*}
where $\G$ is the same set defined in \eqref{eq:gamma}.
\end{lemma}
\begin{lemma}\label{le:>d}
Let $u_{n} \in \H$, $n\ge 1,$ such that $\|u_{n}\|\ge C>0$ and
\begin{equation*}
\max_{t\ge 0}I(tu_{n}) \le c_V+\d_n,
\end{equation*}
with $\d_n\to0^+.$ Then, there exist a sequence $(y_n)_n \subset
\RN$ and two positive numbers $R,\;\mu
>0$ such that
\begin{equation*}
\liminf_n \int_{B_R(y_n)}|u_n|^2 \, d x >\mu.
\end{equation*}
\end{lemma}

\begin{lemma}\label{le:t-bdd}
Let $(u_n)_n  \subset \H$ such that $\|u_n\|=1$ and
\begin{equation*}
I(t(u_n)u_n)= \max_{t\ge 0}I( t u_n) \to c_V,\qquad \hbox{ as }n \to
\infty.
\end{equation*}
Then the sequence $(t(u_n))_n\subset \R_+$ possesses a bounded
subsequence in $\R$.
\end{lemma}

\begin{proof}
We have
\begin{align*}
C\ge \irt |\n u_n|^2 + V(x)u_n^2
=t_n^2 \left(t_n^{p-3}\irt |u_n|^{p+1}- \irt \phi_{u_n}u_n^2\right).
\end{align*}
The conclusion follows from $i$ of Lemma \ref{le:prop} and Lemma
\ref{le:>d}.
\end{proof}

\begin{lemma}\label{le:vn}
Suppose that $V,$ $V_n$ satisfy $(${\bf V1-2}$)$, for all
$n\ge~\!\!1$.
\\
If $V_n\to V$ in $L^\infty(\RN)$ then $c_{V_n}\to c_V.$
\end{lemma}

Now define
\begin{align*}
I_\infty(u) &:=\frac 12 \irt |\n u|^2 + V_\infty u^2
+\frac 14 \irt \phi_u u^2 -\frac{1}{p+1}\irt |u|^{p+1},
\\
c_\infty    &:=c_{V_\infty}.
\end{align*}

As in \cite{R}, we have
\begin{lemma}\label{le:cinfty}
If $V$ satisfies ({\bf V1-3}), we get $c_V<c_\infty$.
\end{lemma}

\begin{proof}
By Theorem \ref{main}, there exists $(w,\phi_w)\in\H\times \D$ a ground state
solution of the problem
\begin{equation*}
\left\{
\begin{array}{ll}
-\Delta u + V_\infty u + \phi u = |u|^{p-1} u & \hbox{in }\RT,
\\
-\Delta \phi = u^2 & \hbox{in }\RT.
\end{array}
\right.
\end{equation*}
Let $t(w)>0$ be such that $t(w)w\in \Ne$. By ({\bf V3}), we have
\begin{align*}
c_\infty & =I_\infty(w) \ge I_\infty \big(t(w) w\big)
\\
&=I\big(t(w) w\big) + \irn \big(V_\infty -V(x)\big) |t(w)w|^2
> c_V,
\end{align*}
and then we conclude.
\end{proof}

\subsubsection{Proof of Theorem \ref{main2}}

Let $(u_n)_n$ be a sequence in $\Ne$ such that
\begin{equation}\label{eq:limN}
\lim_n I(u_n)=c_V.
\end{equation}
Observe that for any $u\in \Ne$ we have
\[
I(u)= \left(\frac 12- \frac{1}{p+1} \right)\irt |\n u|^2+ V(x)u^2
+\left(\frac 14- \frac{1}{p+1} \right) \irt \phi_{u}u^2,
\]
hence, by \eqref{eq:limN}, we deduce that $(u_n)_n$ is bounded in $\H$. Let $\bar u\in\H$ be such that, up to a subsequence,
\begin{align}
&u_n\rightharpoonup \bar u\quad\hbox{weakly in }\H, \label{eq:weak}
\\
&u_n\to \bar u\quad\hbox{in }L^s(B), \hbox{with }B\subset \RT, \hbox{bounded, and }1\le s<6.  \label{eq:loc}
\end{align}
We define the measures
\begin{equation*}
\mu_n(\O)=
\int_{\O} \left(\frac 12- \frac{1}{p+1} \right)
\Big[|\n u_n|^2+ V(x)u_n^2\Big]
+\left(\frac 14- \frac{1}{p+1} \right)
\phi_{u_n}u_n^2.
\end{equation*}
Proceeding as in Lemma \ref{le:concentr}, we infer that
compactness holds for $(\mu_n)_n$, namely there exists a sequence
$(\xi_n)_n$ in $\RT$ with the following property: for any $\d>0$,
there exists $r=r(\d)>0$ such that
\begin{equation}\label{eq:bc3}
\int_{B_r^c(\xi_n)} \left(\frac 12- \frac{1}{p+1} \right)
\Big[|\n u_n|^2+ V(x)u_n^2\Big]
+\left(\frac 14- \frac{1}{p+1} \right)
\phi_{u_n}u_n^2 < \d.
\end{equation}

\noindent{\sc Claim}: $(\xi_n)_n$ is bounded in $\RT$.
\\
Suppose by contradiction that, up to a subsequence, $|\xi_n|\to \infty$, as $n\to \infty$.
\\
Fix $\widehat V<V_\infty$ and let $\widehat I$ be the functional defined as $I$ replacing $V$ by $\widehat V$. For any $n\ge1,$ let $z_n=u_n(\cdot - \xi_n)$ and $\hat t_n>0$ such that the functions $\hat t_n z_n$
are in the Nehari manifold of $\widehat I$.
\\
Let $\d>0$ and consider $r>0$ such that \eqref{eq:bc3} holds. For $n$ sufficiently large, we have
\[
V(x+\xi_n)-\widehat V \ge 0, \qquad \hbox{for all }x\in B_r.
\]
Hence we have
\begin{align*}
c_V+o_n(1) &= I(u_n) \ge I(\hat t_n u_n) = \hat I(\hat t_n u_n) +
\frac{\hat t_n^2}{2}\irt \left(V(x)-\widehat V\right)u_n^2
\\
&\ge c_{\widehat V} +
\frac{\hat t_n^2}{2}\int_{B_r} \left(V(x+\xi_n)-\widehat V\right)z_n^2
+ \frac{\hat t_n^2}{2}\int_{B_r^c} \left(V(x+\xi_n)-\widehat V\right)z_n^2
\\
&\ge c_{\widehat V} - \frac{\hat t_n^2}{2}\int_{B_r^c}
\left|V(x+\xi_n)-\widehat V\right|z_n^2.
\end{align*}
Since by \eqref{eq:bc3}
\[
\int_{B_r^c} \left|V(x+\xi_n)-\widehat V\right|z_n^2\le
C\d,\quad\hbox{for any }n\ge 1,
\]
and $(\hat  t_n)_n$ is bounded (the proof is the same as in Lemma \ref{le:t-bdd}), we get that $c_V \ge c_{\widehat V} - C\d$. By the
arbitrariness in the choice of $\d>0,$ we have $c_V \ge c_{\widehat V}.$
Using Lemma \ref{le:vn} we conclude that $c_V\ge c_\infty,$ which
contradicts Lemma~\ref{le:cinfty}.
\\
\
\\
So $(\xi_n)_n$ is bounded in $\RT$ and then, by \eqref{eq:bc3}, for any
$\delta>0$ there exists $r>0$ such that
\begin{equation}\label{eq:conc}
\|u_n\|_{H^1(B_r^c)}<\d,\quad \hbox{uniformly for }n\ge 1.
\end{equation}
By \eqref{eq:weak}, \eqref{eq:loc} and \eqref{eq:conc}, we have that, taken $s\in [2,6[$,
for any $\d>0$ there exists $r>0$ such that, for any $n\ge 1$ large enough
\begin{align*}
\|u_n-\bar u\|_{L^s(\RT)} & \le \|u_n-\bar u\|_{L^s(B_r)}+\|u_n-\bar u\|_{L^s(B^c_r)}
\\
& \le \d +
C\left(\|u_n\|_{H^{1}(B_r^c)}+\|\bar u\|_{H^{1}(B_r^c)}\right)\le (1+2C)\d,
\end{align*}
where $C>0$ is the constant of the embedding
$H^{1}(B_r^c)\hookrightarrow L^s(B^c_r).$ We deduce that
\begin{equation}\label{eq:conv}
u_n\to \bar u\hbox{ in }L^s(\RT),\;\hbox{for any }s\in [2,6[.
\end{equation}
Since $\phi$ is continuous from $L^{12/5}(\RT)$ to $\D,$
from \eqref{eq:conv} we deduce that
\begin{align}
\phi_{u_n}  \to \phi_{\bar u} \;\hbox{ in }\D,& \qquad \hbox{as } n\to \infty,  \nonumber\\
\irt\phi_{u_n}u_n^2  \to  \irt\phi_{\bar u}\bar u^2,& \qquad \hbox{as } n\to \infty,  \label{eq:convphi}
\end{align}
and for any $\psi\in C_0^\infty(\RT)$
\begin{equation}\label{eq:convint}
\irt\phi_{u_n}u_n\psi  \to  \irt\phi_{\bar u}\bar u\psi.
\end{equation}
By \eqref{eq:limN}, we can suppose (see \cite{Wi}) that $(u_n)_n$
is a Palais-Smale sequence for $I|_{\Ne}$ and, as a consequence, that it is easy to see that
$(u_n)_n$ is a Palais-Smale sequence for $I$.
By \eqref{eq:weak}, \eqref{eq:conv} and \eqref{eq:convint}, we
conclude that $I'(\bar u)=0.$
\\
Since $(u_n)_n$ is in $\Ne$, by 3 of Lemma \ref{le:N}
$(\|u_n\|_p)_n$ is bounded below by a positive constant. As a
consequence, \eqref{eq:conv} implies that $\bar u\neq 0$ and so
$\bar u\in \Ne$.
\\
Finally, by \eqref{eq:weak}, \eqref{eq:conv} and
\eqref{eq:convphi},
\begin{align*}
c_V\le I(\bar u) \le \liminf I(u_n) = c_V,
\end{align*}
so we can conclude that $(\bar u, \phi_{\bar u})$ is a ground state solution of \eqref{SM}.

\section{The critical case}\label{se:crit}

This section is devoted to the study of the critical case and in particular we will give the proofs of Theorem \ref{nonex}, Theorem \ref{main3} and Theorem \ref{main4}.

\subsection{The nonexistence result}

\begin{proofnonex}
Arguing as in \cite{BL1,DM2}, we can prove that if $(u,\phi)\in \H\times \D$ is a solution of the problem \eqref{SMcr}, then $(u,\phi)$ satisfies the following Pohozaev identity:
\begin{equation}\label{eq:Po}
\irt |\n u|^2
+3\irt V(x)u^2
+\irt (\n V(x)\mid x)u^2
+\frac 52 \irt \phi u^2
=\irt u^6.
\end{equation}
Multiplying the first equation of \eqref{SMcr} by $u$ and integrating, we have
\begin{equation}\label{eq:Po2}
\irt |\n u|^2
+\irt V(x)u^2
+\irt \phi u^2
=\irt u^6;
\end{equation}
on the other hand, multiplying the second equation of \eqref{SMcr}
by $\phi$ and integrating, we have
\begin{equation}\label{eq:Po3}
\irt |\n \phi|^2
=\irt \phi u^2.
\end{equation}
By the combination of \eqref{eq:Po}, \eqref{eq:Po2} and \eqref{eq:Po3}, we infer that
\begin{equation*}
2\irt V(x)u^2
+\irt (\n V(x)\mid x)u^2
+\frac 32 \irt |\n \phi |^2 =0,
\end{equation*}
which, together with ({\bf V2}) and ({\bf V4}), implies that $u=\phi=0$.
\end{proofnonex}

\begin{remark}
In fact, the same nonexistence result would hold even if we
supposed the weaker hypothesis
\[
0< C_4 \le 2 V(x) + (\n V(x)\mid x) \le C_5, \qquad \hbox{ for all
}x\in \RT,
\]
in the place of ({\bf V4}).
\end{remark}

\subsection{The existence results}

As in subsection \ref{se:prelim}, for every $u\in L^{12/5}(\RT)$ we
denote by $\phi_u\in \D$ the unique solution of
\[
-\Delta \phi=u^2,\qquad \hbox{in }\RT.
\]
It can be proved that $(u,\phi)\in H^1(\RT)\times \D$ is a
solution of \eqref{SMc} if and only if $u\in\H$ is a critical
point of the functional $I^*\colon \H\to \R$ defined as
\begin{equation*}
I^*(u):= \frac 12 \irt |\n u|^2 + V(x)u^2 +\frac 14 \irt \phi_u
u^2 - \frac {1} {q+1} \irt |u|^{q+1}-\frac{1}6\irt u^6,
\end{equation*}
and $\phi=\phi_u$. \\
The Nehari manifold of the functional $I^*$, defined as
\begin{equation*}
\Ne^*:=\left\{u\in \H\!\setminus\! \{0\} \Big{|} \irt |\n u|^2 +
V(x)u^2 + \phi_u u^2 -   |u|^{q+1} - u^6 = 0 \right\},
\end{equation*}
satisfies the equivalent of Lemma \ref{le:N} and so it is a natural constraint for $I^*$ and we are looking for critical
points of $I^*$ restricted to $\Ne^*.$

Set
\begin{align}
c_1^* &= \inf_{g\in \G^*} \max_{ t\in [0,1]} I^*(g( t)); \nonumber\\
c_2^* &= \inf_{u\neq 0} \max_{ t \ge 0} I^*( t u);\nonumber
\\
c_3^* &= \inf_{u\in\Ne^*} I^*(u);\nonumber
\end{align}
where
\[
\G^*=\left\{g\in C\big([0,1],\H\big) \mid g(0)=0,\;I^*(g(1))\le 0,
\;g(1)\neq 0\right\}.
\]
It is standard to prove that
\begin{lemma}
The following relations hold
\[
c^*_V:=c^*_1=c^*_2=c^*_3.
\]
\end{lemma}

We denote by $S$ the best constant for the Sobolev embedding $\D \hookrightarrow L^6(\RT)$, namely
\[
S  = \inf_{u\in\mathcal{D}^{1,2}\setminus\{0\}}\frac{\|\n u\|_2^2}{\|u\|_6^2}.
\]

\subsubsection{The constant potential case}

In this section we suppose that $V$ is a positive constant. For simplicity we assume $V\equiv 1$ and we denote $c^*=c^*_V$.

\begin{lemma}\label{le:ccc*}
The following inequality holds
\[
c^* <  \frac 13 S^\frac 32.
\]
\end{lemma}

\begin{proof}
Consider the one parameter Talenti's functions
$u_\eps\in\D$ defined by
\begin{equation*}
u_\eps:=C_\eps\frac{\eps^\frac 14}{(\eps + |x|^2)^\frac 12},
\end{equation*}
where $C_\eps>0$ is a normalizing constant (see \cite{T}). Let
$\varphi$ be a smooth cut off function,
namely $\varphi\in C_0^\infty(\RT)$ and
there exists $R>0$ such that $\varphi|_{B_R}=1$, $0\le \varphi\le
1$ and $\supp \varphi\subset B_{2R}.$ Set
$w_\eps:=u_\eps\varphi$ and $v_\eps=w_\eps/\|w_\eps\|_6.$
Using the estimates obtained in \cite{BN} we
get
\begin{align}
\|\n v_\eps\|_2^2 & = S + O(\eps^\frac 12),\label{eq:nabla}
\\ \nonumber
\\
\noalign{\hbox{and, for any $s\in [2,6[,$}}
\|v_\eps\|_{s}^{s} & =
\left\{
\begin{array}{ll}
O(\eps^\frac {s}{4}), & \hbox{if } s\in [2,3[,
\\
O(\eps^\frac{3}{4}|\log(\eps)|), & \hbox{if } s=3,
\\
O(\eps^\frac {6-s}{4}), & \hbox{if } s\in ]3,6[.
\end{array}
\right.\label{eq:s}
\end{align}
For every $\eps>0$ let $ t_\eps > 0$ such that $ t_\eps
v_\eps\in\Ne^*.$ Obviously $( t_\eps)_{\eps >0}$ is
bounded below by a positive constant; otherwise there should exist
a sequence $(\eps_n)_n$ such that $\lim_n  t_{\eps_n}=0$ and then,
by \eqref{eq:nabla}, Lemma \ref{le:prop} and \eqref{eq:s},
\[
0 < c^* \le \lim_n I^*( t_{\eps_n} v_{\eps_n}) = 0.
\]
\\
{\sc Claim:} For any $\eps > 0$ small enough
$ t_\eps \le \left(\irt|\n v_\eps |^2 + v_\eps^2\right)^{1/4}.$
\\
Let $\g_\eps( t):=I^*( t v_\eps)$ and set $r_\eps:=
\left(\irt |\n v_\eps |^2 + v_\eps^2\right)^{1/4}.$ By
\eqref{eq:nabla} and \eqref{eq:s}, $(r_\eps)_{\eps>0}$ is bounded
below by a positive constant. Since $ t_\eps
v_\eps\in\Ne^*,$ certainly $\g'_\eps( t_\eps)=0.$ On the
other hand, by $i$ of Lemma \ref{le:prop} and \eqref{eq:s}, for
any $\eps$ small enough,
\begin{align*}
\g_\eps'( t) &=  t r_\eps^4- t^5+ t^3\irt\phi_{v_\eps}v_\eps^2-
t^q\|v_\eps\|_{q+1}^{q+1}
\\
&\le  t r_\eps^4- t^5 +C' t^3\|v_\eps\|_\frac{12}5^4-
t^q\|v_\eps\|_{q+1}^{q+1}
\\
&= t r_\eps^4- t^5+ t^3 \left(C' O(\eps)- t^{q-3}
O(\eps^{\frac{5-q} 4})\right),
\end{align*}
where $O(\eps)$ and $O(\eps^{\frac{5-q} 4})$ are nonnegative
functions. We deduce that, for any $\eps>0$ small enough,
$\g_\eps'( t)<0$ in $]r_\eps, +\infty [:$ the claim follows as a
consequence.
\\
\\
Now, since the function
\begin{equation*}
 t\in\R_+\mapsto \frac 12  t^2 r_\eps^4
-\frac 16  t^6
\end{equation*}
is increasing in the interval $[0,r_\eps[,$  by \eqref{eq:nabla}
and $i$ of Lemma \ref{le:prop} we have that
\begin{align*}
I^*( t_\eps v_\eps) &
= \frac{ t_\eps^2}2\irt |\n v_\eps|^2 + v_\eps^2
+ \frac{ t_\eps^4}4 \irt \phi_{v_\eps} v_\eps^2
- \frac{ t_\eps^{q+1}}{q+1} \irt |v_\eps|^{q+1} -
\frac{ t_\eps^6}6
\\
&\le \frac 13\left(\irt |\n v_\eps|^2 +
v_\eps^2\right)^\frac 32 +
C'\frac{ t_\eps^4}4\|v_\eps\|_{\frac {12}5}^4
- \frac{ t_\eps^{q+1}}{q+1} \|v_\eps\|_{q+1}^{q+1}
\\
& = \frac 13 \left( S+O(\eps^{\frac 12})
+ \irt v_\eps^2\right)^\frac 32 + C'\frac{ t_\eps^4}4\|v_\eps\|_{\frac {12}5}^4
- \frac{ t_\eps^{q+1}}{q+1} \|v_\eps\|_{q+1}^{q+1}.
\end{align*}
Using the inequality $(a+b)^\d\le a^\d+\d (a+b)^{\d-1}b$
which holds for any $\d\ge 1$ and $a,b\ge 0$, by \eqref{eq:s} and the previous
chain of inequalities we get
\begin{align}\label{eq:S}
I^*( t_\eps v_\eps) \le \frac 13 S^{\frac32}
+ O(\eps^{\frac12}) + C_1(\eps) O(\eps) -
C_2(\eps)O(\eps^\frac{5-q}{4}),
\end{align}
where $C_1(\eps)$ and $C_2(\eps)$ are in an interval
$[\a,\b]$ with $\a>0.$
Since $q>3$, the conclusion follows from \eqref{eq:S}, for $\eps>0$
small enough.
\end{proof}

\begin{proofmain3}
Let $(u_n)_n \subset \Ne^*$ such that
\begin{equation}\label{eq:limmu}
\lim_n I^*(u_n)=c^*.
\end{equation}
We easily deduce that $(u_n)_n$ is bounded in $\H,$ so there
exists $\bar u\in\H$ such that, up to a subsequence,
\begin{align}
&u_n\rightharpoonup \bar u\quad\hbox{weakly in }\H, \label{eq:weakmu}
\\
&u_n\to \bar u\quad\hbox{in }L^s(B), \hbox{with }B\subset \RT, \hbox{bounded, and }1\le s<6.  \nonumber 
\end{align}
As in the first part of the paper, we use a
concentration-compactness argument on the sequence of positive
measures
\begin{align*}
\mu_n^*(\O)= &\left(\frac 12 - \frac{1}{q+1} \right) \int_\O |\n u_n|^2 + u_n^2
+ \left(\frac 14 - \frac{1}{q+1} \right) \int_\O \phi_{u_n} u_n^2
\\
&+ \left(\frac{1}{q+1} - \frac 1 6 \right) \int_\O  u_n^6.
\end{align*}
We define the functional $J^*\colon \H\to\R$ as:
\begin{align*}
J^*(u)=&\left(\frac 12 - \frac{1}{q+1} \right) \irt |\n u|^2 + u^2
+ \left(\frac 14 - \frac{1}{q+1} \right)\irt  \phi_{u} u^2 
\\
& + \left(\frac 1{q+1} - \frac 1 6 \right) \irt u^6.
\end{align*}
\\
\
\\
{\sc Vanishing does not occur}
\\
Suppose by contradiction, that for all $r>0$
\[
\lim_n \sup_{\xi \in \RT}\int_{B_r(\xi)} d \mu_n^* =0.
\]
By \cite{L2} we deduce that $u_n\to0$ in $L^s(\RT)$ for any $s\in ]2,6[.$
\\
By $i$ of Lemma \ref{le:prop}, since $(u_n)_n \subset \Ne^*$, it follows that
\begin{equation*}
\lim_n  \Bigg[ \irt |\n u_n|^2 
+ u_n^2 - \irt u_n^6\Bigg]=0.
\end{equation*}
By the boundedness of $(u_n)_n$ in $\H$, we infer that there exists $l>0$ such that, up to subsequence, 
\begin{equation*}
l:=\lim_n  \irt |\n u_n|^2 + u_n^2 = \lim_n  \irt u_n^6.
\end{equation*}
We have
\begin{equation}\label{eq:prima}
c^*=\lim_n I^*(u_n) = \frac 1 2 \, l - \frac 1 6 \, l = \frac 1 3 \, l
\end{equation}
and
\begin{equation}\label{eq:seconda}
S\le\frac{\irt |\n u_n|^2 + u_n^2}
{{\left( \irt u_n^6 \right) }^\frac 1 3}
\to l ^ \frac 2 3.
\end{equation}
By \eqref{eq:prima} and \eqref{eq:seconda} we get $c^*=\frac 1 3 l
\ge \frac 1 3 S^\frac 3 2,$ contradicting $2$ of Lemma
\ref{le:ccc*}.
\\
\\
{\sc Dichotomy does not occur}
\\
The proof uses similar argument as those in the proof of Theorem
\ref{main}.
\\
\\
So the measures $\mu_n^*$ concentrate and, in particular, we have
that there exists a sequence $(\xi_n)_n$ in $\RN$ such that for
any $\d>0$ there exists $r=r(\d)>0$ such that
\begin{equation}\label{eq:bc2}
\left(\frac 12- \frac{1}{q+1} \right) \int_{B_r^c(\xi_n)}  |\n u_n|^2+ u_n^2  < \d.
\end{equation}
{From} now on, we only give a sketch of the remaining part
of the proof, since it is similar to that of the subcritical case. We define $v_n:=u_n(\cdot -\xi_n)$. It is easy to see that $(v_n)_n\subset \Ne^*$. From \eqref{eq:bc2} we have that for any $\delta>0$
there exists $r>0$ such that
\begin{equation*}
\|v_n\|_{H^1(B_r^c)}<\d,\quad \hbox{uniformly for }n\ge 1.
\end{equation*}
Hence we deduce
\begin{align}
&v_n \to \bar v\hbox{ in }L^s(\RT),\;\hbox{for any }s\in [2,6[;\label{eq:conv2}
\\
&\phi_{v_n}  \to \phi_{\bar v} \;\hbox{ in }\D;  \nonumber
\\
&\irt\phi_{v_n}v_n^2  \to  \irt\phi_{\bar v}\bar v^2. \label{eq:convphi2}
\end{align}
Moreover, for any $\psi\in C_0^\infty(\RT)$,
\begin{align}
\irt\phi_{v_n}v_n\psi  &\to  \irt\phi_{\bar v}\,\bar v \,\psi, \nonumber
\\
\noalign{\hbox{and, by \eqref{eq:conv2},}}
\irt v_n^5 \psi  &\to  \irt \bar v^5\,\psi.\nonumber
\end{align}
By \eqref{eq:limmu}, we can suppose (see \cite{Wi}) that $(v_n)_n$
is a Palais-Smale sequence for $I^*|_{|\Ne^*},$ and, consequently,
it is a Palais-Smale sequence for $I^*.$ By standard arguments, we
infer that $\bar v\in \Ne^*$.
\\
Finally, since $(v_n)_n$ and $\bar v$ are in $\Ne^*$, we have that
\begin{align*}
I^*(\bar v) & = \frac 13 \irt |\n \bar v|^2 +    \bar v^2  
+\frac 1{12} \irt \phi_{\bar v}\bar v^2
+\left(\frac 16 - \frac 1{q+1}\right)\irt |\bar v|^{q+1},
\\
I^*(v_n)  &=   \frac 13 \!\irt |\n v_n|^2 \!+    v_n^2
+\!\frac 1{12} \!\irt \!\phi_{v_n}v_n^2\!
+\!\left(\frac 16 - \frac 1{q+1}\right)\!\!\irt \!|v_n|^{q+1},
\end{align*}
so, by \eqref{eq:limmu}, \eqref{eq:weakmu}, \eqref{eq:conv2} and \eqref{eq:convphi2},
\[
c^* \le I^*(\bar v) \le \liminf I^*(v_n) = c^*.
\]
We conclude that $(\bar v, \phi_{\bar v})$ is a ground state solution of \eqref{SMc}.
\end{proofmain3}

\subsubsection{The non-constant potential case}

In this section we suppose that $V$ satisfies hypotheses ({\bf V1-3}).

We define the functional $I^*_{\infty}:\H\to\R$ and
the Nehari manifold $\Ne^*_{\infty}$ in the following way
\begin{equation*}
I^*_{\infty}(u):= \frac 12 \irt |\n u|^2 + V_\infty u^2 +\frac 14
\irt \phi_u u^2 - \frac {1} {q+1} \irt |u|^{q+1}
-\frac{1}6\irt u^6,
\end{equation*}
\begin{equation*}
\Ne^*_{\infty}:=\left\{u\in \H \! \setminus \! \{0\} \Big{|}
\irt |\n u|^2 + V_\infty u^2 + \phi_u u^2 - |u|^{q+1} -u^6 = 0
\right\}.
\end{equation*}
We set
\[
c_\infty^*=\inf_{u\in \Ne_\infty^*}I^*_\infty(u).
\]

\begin{lemma}\label{le:cS}
The following inequality holds
\[
c^*_V < \frac 13 S^\frac 32.
\]
\end{lemma}
\begin{proof}
By Theorem \ref{main3}, there exists a ground state solution for \eqref{SMc} whenever $V\equiv V_\infty$; so, arguing as in Lemma \ref{le:cinfty}, we can show that $c^*_V < c^*_\infty$. Therefore, the inequality  follows by Lemma \ref{le:ccc*}. 
\end{proof}

Following \cite{R}, by Lemmas \ref{le:ccc*} and \ref{le:cS} and using a non-vanishing type argument as in the proof of Theorem \ref{main2}, we can show that the corresponding versions of Lemmas \ref{le:>d}, \ref{le:t-bdd} and \ref{le:vn} hold for the functional $I^*$.

\begin{proofmain4}
Let $(u_n)_n \subset \Ne^*$ such that
\begin{equation*}
\lim_n I^*(u_n)=c^*_V.
\end{equation*}
We easily deduce that $(u_n)_n$ is bounded in $\H,$ so there
exists $\bar u\in\H$ such that, up to a subsequence,
\[
u_n\rightharpoonup \bar u\quad\hbox{weakly in }\H.
\]
Let us define the positive measure
\begin{multline*}
\mu_n^*(\O)= \left(\frac 12 - \frac{1}{q+1} \right)\int_\O |\n u_n|^2 + V(x)u_n^2
\\
+ \left(\frac 14 - \frac{1}{q+1} \right)\int_\O \phi_{u_n} u_n^2+
\left(\frac{1}{q+1} - \frac 1 6 \right)\irt u_n^6.
\end{multline*}
Arguing as in the proof of Theorem \ref{main3}, we can prove that the measures $\mu_n^*$ concentrate and, in particular, we have
that there exists a sequence $(\xi_n)_n$ in $\RN$ such that for
any $\d>0$ there exists $r=r(\d)>0$ such that
\begin{equation}\label{eq:bc2V}
\left(\frac 12- \frac{1}{q+1} \right) \int_{B_r^c(\xi_n)}  |\n
u_n|^2+V(x)u_n^2  < \d.
\end{equation}
With the arguments similar to those of the subcritical case, we show that the sequence $(\xi_n)_n$ is bounded in $\RT$ and then by \eqref{eq:bc2V} we have that for any $\delta>0$
there exists $r>0$ such that
\begin{equation*}
\|u_n\|_{H^1(B_r^c)}<\d,\quad \hbox{uniformly for }n\ge 1.
\end{equation*}
Hence we deduce
\[
u_n \to \bar u\hbox{ in }L^s(\RT),\;\hbox{for any }s\in [2,6[.
\]
Then, arguing as in the proof of Theorem \ref{main3}, we get the conclusion.
\end{proofmain4}

\end{document}